\documentclass[12pt]{article}
\usepackage{amsmath}
\usepackage{amssymb}
\usepackage{amsfonts}
\usepackage{amsthm}
\theoremstyle{plain}
\newtheorem{theorem}{Theorem}[section]

\newtheorem{proposition}[theorem]{Proposition}

\newtheorem{lemma}[theorem]{Lemma}

\theoremstyle{definition}
\newtheorem{definition}[theorem]{Definition}
\theoremstyle{remark}
\newtheorem*{remark}{Remark}

\newcommand{\refE}[1]{(\ref{E:#1})}
\newcommand{\refS}[1]{Section~\ref{S:#1}}

\newcommand{\refD}[1]{Definition~\ref{D:#1}}

\newcommand{\C}{\ensuremath{\mathbb{C}}}

\renewcommand{\P}{\ensuremath{\mathbb{P}}}
\newcommand{\Z}{\ensuremath{\mathbb{Z}}}

\renewcommand{\i}{{\,\mathrm{i}\,}}
\newcommand{\ka}{{k^{\vee}}}
\renewcommand{\l}{\lambda}
\newcommand{\gb}{\overline{\mathfrak{g}}}
\newcommand{\g}{\mathfrak{g}}
\newcommand{\gh}{\widehat{\mathfrak{g}}}
\newcommand{\ghN}{\widehat{\mathfrak{g}_{(N)}}}

\newcommand{\car}{{\mathfrak{h}}}    
\newcommand{\bor}{{\mathfrak{b}}}    
\newcommand{\nil}{{\mathfrak{n}}}    
\newcommand{\vp}{{\varphi}}
\newcommand{\bh}{\widehat{\mathfrak{b}}}  
\newcommand{\bb}{\overline{\mathfrak{b}}}  
\newcommand{\Vh}{\widehat V}
\newcommand{\KZ}{Kniz\-hnik-Zamo\-lod\-chi\-kov}
\newcommand{\TUY}{Tsuchia, Ueno  and Yamada}
\newcommand{\KN} {Kri\-che\-ver-Novi\-kov}
\newcommand{\pN}{\ensuremath{(P_1,P_2,\ldots,P_N)}}
\newcommand{\xN}{\ensuremath{(\xi_1,\xi_2,\ldots,\xi_N)}}
\newcommand{\lN}{\ensuremath{(\lambda_1,\lambda_2,\ldots,\lambda_N)}}
\newcommand{\iN}{\ensuremath{1,\ldots, N}}
\newcommand{\iNf}{\ensuremath{1,\ldots, N,\infty}}
\newcommand{\MgN}{\mathcal{M}_{g,N}} 
\newcommand{\MgNi}{\mathcal{M}_{g,N}^{\infty}} 
\newcommand{\MgNe}{\mathcal{M}_{g,N+1}} 
\newcommand{\mpt}{(M,P_1,P_2,\ldots, P_N,\Pif)} 
\newcommand{\mpp}{(M,P_1,P_2,\ldots, P_N)} 
\newcommand{\Yt}{\widetilde{Y}}
\newcommand{\bt}{\tilde b}
\newcommand{\sinf}{{\widehat{\sigma}}_\infty}
\newcommand{\St}{\widetilde{S}}
\newcommand{\uni}{\mathcal{U}}
\newcommand{\can}{\mathcal{K}}
\newcommand{\Wh}{\widehat{W}}
\newcommand{\Wt}{\widetilde{W}}
\newcommand{\Pif} {P_{\infty}}
\newcommand{\PN}{\ensuremath{\{P_1,P_2,\ldots,P_N\}}}
\newcommand{\PNi}{\ensuremath{\{P_1,P_2,\ldots,P_N,P_\infty\}}}
\newcommand{\Fl}[1][\lambda]{F^{#1}}
\newcommand{\Fln}[1][n]{F_{#1}^\lambda}
\newcommand{\tang}{\mathrm{T}}
\newcommand{\Kl}[1][\lambda]{\can^{#1}}
\newcommand{\A}{\mathcal{A}}
\newcommand{\U}{\mathcal{U}}
\renewcommand{\O}{\mathcal{O}}
\newcommand{\Ae}{\widehat{\mathcal{A}}}
\newcommand{\Ah}{\widehat{\mathcal{A}}}
\newcommand{\La}{\mathcal{L}}
\newcommand{\Le}{\widehat{\mathcal{L}}}
\newcommand{\Lh}{\widehat{\mathcal{L}}}
\newcommand{\Da}{\mathcal{D}}
\newcommand{\kndual}[2]{\langle #1,#2\rangle}
\newcommand{\cint}{\frac 1{2\pi\i}\int_{C}}
\newcommand{\cintl}{\frac 1{24\pi\i}\int_{C}}
\newcommand{\w}{\omega}
\newcommand{\ord}{\operatorname{ord}}
\newcommand{\res}{\operatorname{res}}
\newcommand{\nord}[1]{:\mkern-5mu{#1}\mkern-5mu:}
\normalbaselines 
\begin{document} 
\vspace*{-1cm}
\hspace*{\fill} Mannheimer Manuskripte 

\hspace*{\fill} math.QA/0001040

\begin{center}
\vspace*{0.9cm}

{\LARGE{\bf 
Elements of a Global Operator Approach to 
Wess-Zumino-Novikov-Witten Models}}

\vskip 1.5cm

{\large {\bf Martin Schlichenmaier }} 

\vskip 0.4 cm 

Department of Mathematics and Computer Science\\ 
University of Mannheim, D7, 27 \\ 
D-68131 Mannheim, Germany

\end{center}

\vspace{0.3 cm}

\begin{abstract}
Elements of a global operator approach to the WZNW models
for compact Riemann surfaces of arbitrary genus $g$  with $N$ 
marked points were given by Schlichenmaier and Sheinman.
This contribution reports on the results.
The approach is based on the multi-point Krichever-Novikov
algebras of global meromorphic functions and vector fields, and 
the global algebras of affine type 
and their representations.
Using the global Sugawara construction and the
identification of a certain subspace of the vector field algebra
with the tangent space to the moduli space of the geometric data,
Knizhnik-Zamalodchikov equations are defined. Some steps
of the approach of Tsuchia, Ueno and Yamada 
to WZNW models are presented to compare
it with our approach.
\end{abstract}

\vspace*{0.2cm}
{\sl  Invited talk  presented at the 
3rd International Workshop on
"Lie Theory and Its Applications in Physics - Lie III",
11 - 14 July 1999, Clausthal, Germany.}

\vspace{0.5 cm} 

\section{Introduction}
\label{S:intro}
Wess-Zumino-Novikov-Witten (WZNW) models provide important examples of 
a two-dimensional conformal field theory.
They can roughly be described as follows. The gauge algebra of the theory is
the affine algebra associated to a finite-dimensional gauge
algebra (i.e. a simple finite-dimensional Lie algebra). The
geometric data consists of a compact Riemann surface (with complex 
structure) of genus $g$ and a finite number of marked points on 
this surface.
Starting from representations 
of the gauge algebra the space of conformal
blocks can be defined. 
It depends on the geometric data.
Varying the geometric data should yield a bundle 
over the  moduli space of the geometric data. 
In \cite{KniZam} Knizhnik and Zamolodchikov 
considered the case of genus 0 (i.e. the Riemann sphere).
 There, changing the geometric data consists in moving 
the marked points on the sphere. The space of conformal blocks 
could completely be found inside the part of the representation associated
to the finite-dimensional gauge algebra.
On this space an important set of equations, 
the Knizhnik-Zamolodchikov (KZ) equations, was introduced.
In a general geometric setting, solutions are the  flat sections of
the bundle of conformal blocks over the moduli space with respect
to the KZ connection.

For higher genus it is not possible to realize the space of
conformal blocks inside the representation space associated to
 the finite-dimensional
algebra. There exists different attacks to the
generalization. Some of them add additional structure on
these representation spaces (e.g. twists, representations of
the fundamental group,..).
Here I do not have the place to pay  proper reference to all these 
approaches. Let me only give a few names:
Bernard \cite{Ber1},\cite{Ber2}, Felder and Wieczerkowski \cite{Fel},
\cite{FeWi},  Hitchin \cite{Hi}, and Ivanov \cite{Iv}.

An important approach very much in the spirit of the original
\KZ\ approach was given by Tsuchia, Ueno and Yamada \cite{TUY}.
In \refS{tuy} I will present a very short outline of their theory.
The main point in their approach is that at the marked points, 
after choosing local coordinates, local
constructions  are done.
In this setting
the well-developed theory of representations of the traditional
affine Lie algebras (Kac-Moody algebras of affine type) can be used. 
It appears a mixture between local and global objects and
considerable effort is necessary to extend the local constructions
to global ones.

Oleg Sheinman and myself propose a different approach to the
WZNW models which uses consequently only global 
objects.
These objects are the Krichever-Novikov (KN) algebras \cite{KN} 
and their representations, respectively  their multi-point generalizations
given by me \cite{Schlce}. An outline of this approach is presented in 
\refS{goa}.
In the remaining sections more details on the construction are given.
Here I only want to point out 
that a subspace of the KN algebra of vector fields is identified
with tangent directions on the moduli space of the geometric data.
Conformal blocks can be defined.
We are able to incorporate a richer theory because in our
set-up we are able to deal with more general representations of the
global algebras.

Finally, with the help of the global Sugawara construction for 
the higher genus and
multi-point situation proven in \cite{SchlShSug}, it is possible to
define the higher genus multi-point Knizhnik-Zamolodchikov equations
(see \refD{KZe}).
The complete proofs of the results appeared in \cite{SchlShWZ1}.
A detailed study of the KZ equations, resp. of the connection 
is work in progress \cite{SchlShWZ2}.
\section{Outline of the Tsuchya-Ueno-Yamada approach} 
\label{S:tuy}
Let me first recall a few steps of the approach of \TUY\ to the WZNW models.
I will concentrate on the steps which are of relevance in our approach.
More details can be found in \cite{TUY}, resp. in the more pedagogical
introduction \cite{Ucft}.

{\bf (1)} 
A finite-dimensional complex simple Lie algebra 
$\g$ 
with $(.|.)$ a symmetric, nondegenerate, invariant bilinear form
(e.g. the Killing form) is fixed.
 This Lie algebra is the finite-dimensional gauge algebra of the theory.
Some authors prefer to call this algebra the horizontal algebra.

{\bf (2)} 
Fix a natural number $N$ and  a compact Riemann surface $M$ of arbitrary
genus $g$ (resp. a smooth projective curve over $\C$).
Take a tuple of $N$ distinct points 
$\pN$ on $M$ and around every such point $P_i$ a coordinate $\xi_i$.
This defines the geometric data 
\begin{equation}
  \label{E:gdcoor}
  \Upsilon'=(M,\pN,\xN)\ .
\end{equation}

{\bf (3)}
The affine algebra (or Kac-Moody algebra of affine type) 
associated to $\g$ is given as
\begin{equation}
  \gh=\g\otimes\C\,[t^{-1},t]]\oplus \C\, c
\end{equation}
with Lie structure
\begin{gather}
[x\otimes t^n,y\otimes t^m]=[x,y]\otimes t^{n+m}+
(x|y)\, n\cdot \delta_{m}^{-n}\, c,
\qquad \forall n,m\in\Z\ ,
\\
[x\otimes t^n,c]=0,\qquad 
\forall n\in\Z\ .  
\end{gather}
Note that in this approach one has to consider Laurent series 
instead of the usual Laurent polynomials. In our approach 
we will return to Laurent polynomials.
By ignoring the central element $c$ one sees that
$\gb=\g\otimes \C\,[t^{-1},t]]$ carries also a Lie algebra structure.
This algebra is called the loop algebra (associated to $\g$).
There is a well-defined theory of highest weight representations
$H_\l$ 
of the affine algebra 
$\gh$
associated to a level and certain weights $\l$ of the finite
dimensional Lie algebra $\g$, see \cite{KacB}.
Recall that the central element $c$ operates as \ level$\times id$\  on 
$H_\l$.

Another Lie algebra appearing in this context is the
Virasoro algebra $V$ which is the Lie algebra with basis 
$\{l_n,n\in\Z\}\cup \{c_1\}$ and  
Lie structure
\begin{gather}
[l_n,l_m]=(m-n)\, l_{n+m}+\delta_{m}^{-n}\frac {n^3-n}{12}\,c_1,
\qquad \forall n,m\in\Z\ ,
\\
[l_n,c_1]=0,\qquad 
\forall n\in\Z\ .  
\end{gather}
Starting from a highest weight representation of the affine algebra
the Sugawara construction 
defines also a representation of the  Virasoro algebra on the
representation space.
Further down I will describe the Sugawara construction 
in a general geometric setting.
This setting will incorporate also the classical Sugawara construction.
Note that in the  theory a graded structure is 
implicitly given and is employed.
For example the Virasoro algebra becomes a graded affine Lie algebra
by defining $\deg(l_n):=n$ and $\deg(c_1):=0$.

{\bf(4)}
After fixing the level globally,
one assigns to every marked point $P_i$ a 
heighest weight representation
$H_{\l_i}$ of the gauge algebra.
Set $\l:=\lN$ and 
\begin{equation}
H_\l:=H_{\l_1}\otimes H_{\l_2}\cdots\otimes H_{\l_N}\ .
\end{equation}
The space $H_\l$ is a representation space for the Lie algebra $\ghN $,
where this algebra is defined as the one-dimensional central
extension of $N$ copies of the loop algebra. 
Here the $i$-th copy of the loop algebra operates on the
$i$-th 
factor in the tensor product and is associated to the $i$-th
point $P_i$.
To distinguish the different copies we use 
$t_i$ for the affine parameter $t$ in the loop algebra.

{\bf(5)}
Up to now the geometry was not really involved. 
This changes in this step.
The affine parameter $t_i$ corresponding formally to  $P_i$ via
the assignment of the representation $H_{\l_i}$ to this point,
will be identified with the coordinate $\xi_i$ at this point.
One considers the 
{algebra} $A(\Upsilon')$ of meromorphic functions on $M$ 
which have poles at most at the points $\{P_1,P_2,\ldots,P_N\}$,
and sets
$\g(\Upsilon'):=\g\otimes A(\Upsilon')$.
This algebra is also called the block algebra.
By taking the Laurent expansion of $f\in A(\Upsilon')$ at the point $P_i$ 
with respect to the coordinate $\xi_i$ there, we get an embedding
$  \g(\Upsilon')\to \gb_N\ $
by assigning  to $f$ the corresponding
Laurent series in the affine parameter $t_i$ in the $i$-th copy of $\gb$.
The cocycle defining the central extension vanishes on 
$\g(\Upsilon')$. Hence it can be considered as  
a subalgebra of $\ghN $.

{\bf(6)}
The space of conformal blocks (also called chiral blocks) are defined as the
coinvariants
\begin{equation}
  \label{E:clcb}
  V_{\l}:=H_{\l}/\g(\Upsilon')H_{\l}\ .
\end{equation}
In some context is is better to work with the dual
objects
$   V_{\l}^{*}=\mathrm{Hom}_{\C}(V_\l,\C)$.
This space can be described as the space of linear forms on 
$H_{\l}$ vanishing on $\g(\Upsilon')H_{\l}$.
The vector spaces $ V_{\l}$ turn out to be  finite-dimensional. Their
dimension is given by the Verlinde formula.

{\bf(7)}
One of the  motivations of \cite{TUY} was to 
supply a proof of the Verlinde formula.
For this the authors  pass to the moduli space $\MgNi$ of the 
data $\Upsilon'$ (with the obvious identifications under isomorphisms).
Clearly, not only the goal to proof the Verlinde formula forces
one to consider moduli spaces but also the concept of
quantization requires to consider all possible configurations.
Now everything has to be sheavified. One obtains the
sheaf of conformal blocks over moduli space. On this sheaf the
\KZ\ connection is constructed. The sheaf is indeed
a vector bundle, hence the dimension of the conformal blocks
will be constant along the moduli space,
only depending on the genus $g$, the number  $N$ of 
marked points, and 
the associated weights $\l$  (of course the Lie algebra 
$\g$ and the level will be fixed).
The construction of the connection involves the Sugawara
construction which is done locally on the Riemann surface.
By introducing a projective connection on 
the Riemann surface it is possible to show that it 
globalizes.

{\bf(8)}
It is shown that the essential data can be given in terms of the moduli space
$\MgN$ 
of smooth projective curves of genus $g$ with $N$ marked points,
by forgetting the coordinate systems at the marked points.

{\bf(9)}
Finally, by passing to the Deligne-Mumford-Knudsen boundary of
$\MgN$ corresponding to stable singular curves
and by proving factorization rules (e.g. behaviour of the 
conformal blocks under normalization of the singular curves)
it is possible to express the dimension of the space of conformal blocks 
for the $(g,N,\l)$ situation in dimensions of  spaces of conformal
blocks for lower genera.
\section{Outline of the global operator approach} 
\label{S:goa}
The approach presented in \refS{tuy}
is very successful and a beautiful piece of mathematics.
Nevertheless it has some problems.

{\bf(a)}
It is necessary to choose coordinates around
the points $P_i$. The authors have to work over 
infinite-dimensional moduli spaces. Finally, a complicated reduction
process is needed to reduce the relevant
data to the moduli space of  curves with marked points.
Note that by Tsuchimoto \cite{Tsuch} some simplifications with respect to
this problem has been given.

{\bf(b)}
The proof that the local Sugawara construction globalizes 
(with the help of a projective connection) 
on the Riemann surface is difficult.

{\bf(c)}
Some objects of the theory are only defined locally, other
objects globally.

{\bf(d)}
The representation  $H_\l$ of the affine algebra $\ghN$ does not
see the geometry.

\medskip
We (Oleg Sheinman and myself) propose a different approach
to the  WZNW models which still stays in the algebraic-geometric
setup.
In this section I will outline what has been done so far. Some more details
will be given in the following sections.
First note that the geometric data we consider consist 
of a compact Riemann surface $M$ (or a smooth projective curve over 
$\C$ -- I will use the terms interchangeable) 
and $N$ (distinct) marked points:
\begin{equation}
  \label{E:gd}
  \Upsilon=(M,\pN)\ .
\end{equation}
In contrast to the data \refE{gdcoor}
it does not contain coordinates.
Again, we fix a simple finite-dimensional complex Lie
algebra $\g$.

{\bf(1)}
We replace all local algebras by algebras of \KN (KN)-type
and their multi-point generalizations.
They consists of meromorphic objects
which might have  poles only at the points $\{P_1,P_2,\ldots,P_N\}$ and
a fixed  reference point $\Pif$.
This can be done for the function algebra, 
the affine algebra, the vector field algebra,
the differential operator algebra.
Central extensions can be defined.
They carry an almost-graded structure (see the definition below)
induced by the vanishing  order at $\{P_1,P_2,\ldots,P_N\}$, resp.
at $\Pif$.
In particular, all these algebras can be decomposed 
into subalgebras of elements of positive degree 
 corresponding roughly to elements holomorphic
at $\{P_1,P_2,\ldots,P_N\}$ and subalgebras of elements of negative degree
corresponding roughly to elements holomorphic at $\Pif$ and a
``critical'' finite-dimensional subspace  spanned by other elements.

{\bf(2)}
In the next step we consider highest weight representations of the
higher genus affine algebras. Examples 
are given again by representations $\Wh_\l$
which are constructed from  weights assigned to the marked points.
For highest weight  representations the Sugawara construction works and 
gives a representation of the centrally extended vector field algebra
(see \cite{SchlShSug}, \cite{SchlSug}).

{\bf(3)}
The objects make direct sense over a dense subset of $\MgN$.
We obtain sheaf versions of our objects.

{\bf(4)}
Note that the block algebra $\g(\Upsilon)$ considered above is naturally
a subalgebra of our higher genus multi-point affine algebra.
There is no need to choose coordinates.
It is possible to define again 
\begin{equation}
  \label{E:cbn}
  \Vh_{\l}:=\Wh_{\l}/\g(\Upsilon)\Wh_{\l}\ .
\end{equation}
as the space of conformal blocks.

{\bf(5)}
At a fixed point $b$ of the moduli space $\MgN$ 
the basis elements $l_k$ of the 
``critical subspace''
of the vector field algebra
can be identified with the tangent directions $X_k$ to  
the moduli space at $b$. Additionally the elements operate
as Sugawara operators $T[l_k]$ on the 
representations $\Wh_\l$ 
defined above $b$.
This allows us to define the formal \KZ (KZ) equations to be 
\begin{equation}
  \label{E:kz}
  (\partial_k +T[l_k])\,\varPhi=0,\qquad k=1,\ldots, 3g-3+N\ .
\end{equation}
Here $\varPhi$ is assumed to be a section of $\Wh_\l$ , resp. of
$\Vh_\l$, and $\partial_k$ is an action of $X_k$ on $\Wh_\l$
which fulfills  suitable conditions.

In the following sections I will give additional information 
on these steps. Details appeared in 
\cite{SchlShWZ1}.
Further work is in progress
\cite{SchlShWZ2}.
In particular what has to be done is to construct a 
projectively flat connection on the sheaf of conformal blocks and to
extend the construction to the whole moduli space. Especially it should
be extended its boundary to obtain again a proof of the Verlinde formula.
\section{The Krichever-Novikov objects} 
\label{S:kno}
Let $M$ be a compact Riemann surface of genus $g$.
Let $A$ be a fixed set of finitely many points on $M$ which 
is splitted into two non-empty subsets $I$ and $O$.
This is the general situation dealt with in 
\cite{Schlce}. Here it is enough to consider
\begin{equation}
  \label{E:split}
  I:=\PN,\qquad O:=\{\Pif\}\ .
\end{equation}
The classical (Virasoro) situation is: 
$\   M=\P^1\ $  with $I=\{z=0\}$ and $O=\{z=\infty\}$.
The \KN\ situation \cite{KN}  is:
$  M $ an arbitrary Riemann surface  with  $I=\{P_+\}$ and
$O=\{P_-\}$.

Let $\can$ be the canonical bundle,
i.e. the bundle whose local sections are the holomorphic 
1-differentials, 
 and $\Kl$ (for $\l\in\Z$) its tensor
powers,
resp. for $\l$ negative, the tensor power of its dual.
 Denote by $\Fl(A)$  the space of meromorphic sections of $\Kl$
consisting of those elements which are holomorphic outside of $A$.
We set $\A(A):=\Fl[0](A)$ for the associative 
algebra of functions and $\La(A):=\Fl[-1](A)$
for the Lie algebra of vector fields
(with the usual Lie bracket as Lie structure).
If the set $A$ is clear from the context, we will drop it in the notation.
By multiplication with the elements of $\A$ the space $\Fl$ 
becomes a module over $\A$. By taking the Lie derivative with respect
to the vector fields it becomes also  a Lie module over $\La$.
\begin{theorem}
  \label{T:decomp}
(\cite{Schlce}, \cite{Schlth})
There exists a decomposition 
\begin{equation}
  \label{E:decomp}
  \Fl=\bigoplus_{n\in\Z} \Fln\ ,
\end{equation}
where the $\Fln$ are subspaces of dimension $N$.
By defining the elements of $\Fln$ to be the homogeneous 
elements of degree $n$, the algebras  
$\A$ and $\La$ are almost-graded (Lie) algebras
and the $\Fln$ are almost-graded (Lie) modules over 
$\A$, respectively over $\La$.
\end{theorem}
\noindent
For the convenience of the reader let me recall the
definition of an almost-grading at the
example of the Lie algebra $\La$.
A Lie algebra $\La$ is called almost-graded if it can be decomposed
(as vector space)
\begin{equation}
 \La=\bigoplus_{n\in\Z} \La_n,\quad \dim \La_n<\infty\ ,
\end{equation}
such that there exists $L,K\in\Z$ with
\begin{equation}
\label{E:almgr}
  [\La_n,\La_m]\quad\subseteq \quad\bigoplus_{h=n+m-L}^{n+m+K} \La_h,\quad
\forall n,m\in\Z\ .
\end{equation}
This definition (suitable modified) works for the associative algebra 
$\A$ and the modules $\Fl$.
In our situation we 
always can do with lower shift $L=0$.

The almost-grading is important to obtain a decomposition of the 
algebras into positive and negative parts and to define the
concept of highest weight representations (see below).
It is also necessary to obtain an embedding
of the algebras $\A$ and $\La$ (and more generally also of $\Da$,
the Lie algebra of differential operators, see \cite{Schlct1})
into $\overline{gl}(\infty)$ via their action on the modules $\Fl$.
Here $\overline{gl}(\infty)$ denotes the algebra of (both-sided) 
infinite matrices with finitely many diagonals.

The degree is introduced by the order of $f\in\Fl$ at $I$.
More precisely, $\Fl_n$ is given by exhibiting a basis 
$f_{n,p}^\l,\ p=\iN$. For the generic situation ($\l\ne 0,1$, 
$g\ne 1$, and a generic choice of $A$) 
the element is fixed up to multiplication with a scalar
by
\begin{equation}
\ord_{P_i}(f_{n,p}^\l)=
\begin{cases}
  n-\l,&i=p
\\
 n-\l+1,&i\ne p, i=1,\dots,N
\\
-N(n+1-\l)+(2\l-1)(g-1),&i=\infty\ .
\end{cases}
\end{equation}
For the non-generic situation for finitely 
many $n$ the prescription at $i=\infty$ 
has to be adjusted.
A detailed prescription is given in 
\cite{Schlce} and \cite{Schlth}.
To fix the element 
$f_{n,p}^\l$ uniquely it is necessary to fix a coordinate $\xi_i$ 
centered at 
$P_i$,  for $i=\iN$ and to require
$f_{n,p}^\l(\xi_p)_{|}=\xi_p^{n-\l}(1+O(\xi_p))d\xi_p^\l$.
Note that the fixing does not really depend on the full information
in the coordinate. It depends only on the first infinitesimal 
neighbourhood.
In detail: Let 
$\xi_p'=a_1\xi_p+\sum_{j\ge 2}a_j\xi_p^j$
(clearly with $a_1\ne 0$) 
be another coordinate centered at $P_p$ then 
the normalization  will only depend on the 
value $a_1$.

We have the important duality
\begin{equation}
  \label{E:kndual}
  \Fl\times \Fl[{1-\l}]\to \C,\quad
\kndual f g=\cint f\cdot g=-\res_{\Pif}(f\cdot g),
\end{equation}
where $C$ is any separating cycle for $\pN$ and $\Pif$ which is
cohomologous to (-1) $\times$ the circle around $\Pif$.
For the scaled basis elements we obtain
\begin{equation}
  \label{E:knd}
  \kndual {f_{n,p}^\l}{f_{m,r}^{1-\l}}=
\delta_n^{-m}\delta_p^r\ .
\end{equation}
Again note that the duality relation \refE{knd}
will not depend on the coordinates
chosen.
For special values of $\l$ we set
\begin{equation}
  A_{n,p}:=f_{n,p}^0,\quad e_{n,p}:=f_{n,p}^{-1},\quad
 \w^{n,p}:=f_{-n,p}^{1},\quad \Omega^{n,p}:=f_{-n,p}^{2} .
\end{equation}
In particular $(A_{n,p},\w^{n,p})$ and 
$(e_{n,p},\Omega^{n,p})$ are dual systems of basis elements.

Next we decompose our algebras.
First we consider $\A$, the algebra of functions:
\begin{gather*}
  \A=\A_-'\oplus \A_{(0)}'\oplus\A_+,
\\
\A_-'=\bigoplus_{n\le-K-1}\A_n,\quad
\A_+=\bigoplus_{n\ge 1}\A_n,\quad
\A_{(0)}'=\langle A_{n,p}\mid 1\le p\le N,\ -K\le n \le 0\rangle
\ .
\end{gather*}
Here  $\A_n$ is the subspace of $\A$ consisting of the elements of
degree $n$ and 
$K$ is the constant appearing in the definition of
the almost-grading \refE{almgr} for the algebra $\A$.
The constant depends on the genus $g$ and the number of points $N$.
{}From the almost-grading it follows that $\A_+$ and 
$\A_-'$ are subalgebras. In general $\A_{(0)}'$ is only a subspace.
For our purpose here it is more convenient 
to take as $\A_-\supseteq\A_-' $
  the subalgebra of meromorphic functions vanishing at
$\Pif$,
make a change of basis 
in the lowest degree part $\A_{-K}$ in  $\A_{(0)}$, and take 
$\A_{(0)}\subseteq \A_{(0)}'$ correspondingly smaller.
This yields the decomposition
\begin{equation}
  \label{E:adecomp}
  \A=\A_-\oplus \A_{(0)}\oplus\A_+\ .
\end{equation}
A completely analogous decomposition exists for the vector field
algebra $\La$
\begin{equation}
  \label{E:ldecomp}
  \La=\La_-\oplus \La_{(0)}\oplus\La_+\ .
\end{equation}
The vector fields  in $\La_+$ are vanishing of order 2 at the points  $\PN$,
and the vector fields in   $\La_-$ are vanishing of order 2 at the point
$\Pif$.
The space $\La_{(0)}$ is $(3g-3 +2N +2)$-dimensional. We call it the
``critical strip''.

For further reference let me identify the basis elements of the
critical strip.
First we have $N$ elements $e_{0,p},\ p=\iN$ vanishing at all $P_i,\ i=\iN$,
but not vanishing of 2nd order at every point.
Next we have $N$ elements $e_{-1,p},\ p=\iN$ regular at all $P_i,\ i=\iN$,
but not vanishing   at every point.
Later we will see that they will be responsible for moving the points
$P_p,\ p=\iN$.
In the middle  we have $3g-3$ elements which are neither regular at
$\PN$ nor at $\Pif$.
They will correspond to  deformations of the complex structure of $M$.
Finally we  have two elements regular of order 1, resp. of order 0 at
$\Pif$.

For the algebras we can construct central extensions
$\Ae$ and $\Le_R$ defined via the following geometric
cocycles
\begin{align}
\label{E:centf}
\A:& \qquad
\gamma(g,h):=\cint gdh,
\\
\label{E:centv}
\La:&\qquad
\chi_R(E,f):=\cintl
\left(\frac 12\left(e'''f-ef'''\right)-R\left(e'f-ef'\right)\right)dz\ .  
\end{align}
Here $R$ is a projective connection on $M$.
The cocycles are local in the sense that there exists $T$ and $S$ such
that for all $n,m\in\Z$: 
$\gamma(\A_n,\A_m)\ne 0\implies 
T\le m+n\le 0$ and 
$\chi_R(\La_n,\La_m)\ne 0\implies 
S\le m+n\le 0$.
Again the $T$ and $S$ can be explicitly calculated \cite{Schlth}.
Restricted to the $(+)$ and $(-)$ subalgebras in 
\refE{adecomp} and \refE{ldecomp}
the cocycles are vanishing.
In particular, the subalgebras can be considered in a natural way as
subalgebras of the central extensions $\Ae$ and $\Le_R$.
By setting $\deg(t):=0$ for the central elements the almost-grading
can be extended to the central extensions.

Different connections $R$ 
yield cohomologous cocycles $\chi_R$, hence equivalent
central extensions (see \cite{Schlth} for details).
In the classical situation of $\P^1$ with two marked points one obtains the
Virasoro algebra.
\section{The affine multi-point Krichever-Novikov algebras} 
\label{S:akn}
Let $\g$ be a finite-dimensional reductive Lie algebra
(e.g. semi-simple or abelian) and $(.|.)$  a 
symmetric, nondegenerate invariant bilinear form on $\g$.
The (higher genus multi-point) loop algebra is defined to be 
\begin{equation}
  \label{E:loop}
  \gb:=\g\otimes \A\qquad\quad\text{with}\quad
[x\otimes g,y\otimes h]:=[x,y]\otimes (g\cdot h)\ .
\end{equation}
Its elements  can be considered as $\g$-valued meromorphic functions on
$M$ which are holomorphic outside of $\PNi$.
A central extension $\gh$ is given as 
vector space $\gh=\C\oplus\gb$
with Lie structure  ($\widehat{x}:=(0,x), \ t:=(1,0)$)
\begin{equation}
  \label{E:knaff}
  [\widehat{x\otimes f}, \widehat{y\otimes g}]
 :=\widehat{[x,y]\otimes (f\cdot g)}-(x,|y)\,\gamma(f,g)\cdot t\ ,
\qquad [t,\gh]=0\ ,
\end{equation}
and $\gamma$ from \refE{centf}.
We call this algebra higher genus multi-point 
affine algebra. For $g=0$ and 2 points the classical affine algebras
are obtained (with respect to  Laurent polynomials).
For higher genus and two points it was introduced by \KN\ \cite{KN}
and extensively studied by Sheinman
\cite{Shhwe},\cite{Shhwe2},\cite{Shhwe3}.
Its generalization to higher genus and an arbitrary number of marked
points was given in \cite{Schlct1}.

By the decomposition \refE{adecomp} of $\A$ we obtain a decomposition
\begin{equation}
  \label{E:afdecomp}
  \gh=\gh_-\oplus \gh_{(0)}\oplus\gh_+\ .
\end{equation}
In particular, $\gh_{(0)}=\A_{(0)}\otimes \g\oplus \C \,t$.
Using the duality \refE{knd} we see 
that $1=\sum_{p=1}^N A_{0,p}$ and hence 
that via $\g\cong\g\otimes 1\hookrightarrow \gh$
the finite dimensional Lie algebra 
$\g$ can be naturally considered as subalgebra of 
$\gh_{0}\oplus\gh_+$, of  
$\gb$, and of $\gh$.
\section{Verma modules and Sugawara construction} 
\label{S:verma}
In this section we assume $\g$ to be a simple complex Lie algebra.
Choose $\car\subseteq\g$ a Cartan subalgebra, $\bor$ a 
corresponding Borel subalgebra, and 
$\nil$ a corresponding upper nilpotent subalgebra.
We choose $\l=\lN,\l_p\in\car^*$, take $N$ copies of everything,
and label them by $p=\iN$.
We set
\begin{equation}
  \g_{(N)}:=\g_1\oplus\cdots\oplus\g_N,\quad
  \bor_{(N)}:=\bor_1\oplus\cdots\oplus\bor_N\ .
\end{equation}
Let $W_p$ be a one-dimensional vector space over $\C$ with basis $w_p$
for $p=1,\ldots,N$.
On $W_p$ a one-dimensional representation of $\bor_p$ is defined via
\begin{equation}
  h_p\w_w:=\l_p(h_p)w_p,\quad
  n_pw_p:=0,
  \qquad
h_p\in\car_p,\ n_p\in\nil_p\ .
\end{equation}
This yields a representation of $\bor_{(N)}$
on $W=\bigotimes_{p=1}^N W_p$.
As usual the Lie algebra acts via ``Leibniz rule'' on the
tensor product, and on the  $p$-th tensor factor 
only the $p$-th summand operates non-trivially.
We call this representation $\tau_\l$.
\begin{remark}
It is possible to allow twists to get a richer theory.
They can be obtained via automorphisms $\vp_p:\g\to\g=\g_p$ given
by
$\vp_p=Ad\,\gamma_p,\ \gamma_p\in G$, where $G$ is the Lie group
associated to $\g$.
\end{remark}
\begin{lemma}
The subspace  $\ \gb_0\oplus\C \,t\oplus\gh_+\ $ is a Lie subalgebra 
of $\gh$ and 
\newline
$\psi:\gb_0\oplus\C\, t\oplus\gh_+\to\ghN$ 
defined by
\begin{equation}
  \psi(\sum_{p=1}^N x_p\otimes A_{0,p}):=(x_1,\ldots,x_N),
\quad \psi(t):=0,\quad \psi(\gh_+):=0
\end{equation}
is a Lie homomorphism.
\end{lemma}
\begin{proof}
We calculate
$$
[x_p\otimes A_{0,p},y_r\otimes A_{0,r}]
=[x_p,y_r]\otimes(A_{0,p}A_{0,r})=
[x_p,y_r]\otimes
(A_{0,r}\delta_p^r+\sum_{h>0}\sum_{s=1}^N \alpha_{(0,p),(0,r)}^{(h,s)}
A_{h,s})\ ,
$$
with $\alpha_{(0,p),(0,r)}^{(h,s)}\in\C$.
First, this shows that the subspace is indeed closed under the 
Lie bracket.
Note that  the  defining cocycle for the central
extension vanishes for the subalgebra.
Next, we see that all mixing terms are of higher degree. 
They will be annulated under
$\psi$. Only the term with $p=r$ will survive.
This shows that $\psi$ is a Lie homomorphism.
\end{proof}
Denote by $\bh$ the subalgebra  $\bb_0\oplus\C\, t\oplus\gh_+$
(recall $\bb_0=\bor\otimes \A_0$).
By restricting $\psi$ to $\bh$ we obtain from the lemma that 
$\psi: \bh\to \bor_{(N)}$
is a Lie homomorphism.
Now we  can define 
for $\delta\in\C$ a representation of $\bh$ on
$W$ by
\begin{equation}
  \tau_{\l,\delta}(x\oplus \alpha t\oplus x_+):=
\tau_\l(\psi(x))+\alpha\delta\cdot id\ .
\end{equation}
\begin{definition}
The linear space 
\begin{equation}
  \label{E:verma}
  \widehat{W}_{\l,\delta}:=U(\gh)\otimes_{U(\bh)} W
\end{equation}
with its natural structure of a $\gh$-module is called the Verma module
of the Lie algebra $\gh$ corresponding to 
$(\l,\delta)$ where  $\l\in(\car^*)^N$ 
is the weight of the Verma module and $\delta\in\C$ is the level of
the Verma module.
\end{definition}
\noindent
As usual $U(\gh)$ denotes the universal enveloping algebra.

Now the question arises can the finite-dimensional representations of
$\g$ be recovered. The answer is yes. It is given by
\begin{proposition}
The $\gh$-module $\widehat{W}_{\l,\delta}$ is under the natural
embedding of $\g$ into $\gh$ also a $\g$-module and contains 
the  $\g$-module
\begin{equation}
  \Wt_\l=\Wt_{\l_1}\otimes\cdots \Wt_{\l_N}\ ,
\end{equation}
where the  $\Wt_{\l_i}$
are the heighest weight modules of $\g$ of weight $\l_i$. 
\end{proposition}
The proof is straightforward.
The module  $\Wt_\l$ lies in the degree zero part of $\widehat{W}_{\l,\delta}$,
but there are other elements of degree zero not corresponding to  $\Wt_\l$.

Denote by $\gh_-^*=\gb_-^*\subseteq\gh$ the subalgebra of 
$\g$-valued meromorphic
functions which are regular at $\Pif$,
then we can define the space of conformal blocks as the
space of coinvariants
\begin{equation}
  \label{E:cb}
  \Vh_{\l,\delta}:=\widehat{W}_{\l,\delta}/\gh_-^*\widehat{W}_{\l,\delta}\ .
\end{equation}
Verma modules are examples of admissible modules $\Wh$ 
of $\gh$ in the
following sense:
\begin{enumerate}
\item 
the central element $t$ operates as $c\cdot id$ with $c\in\C$,
($c$ is called the level),
 \item
for all $w\in\Wh$ we have $\gh_n.w=0$,  for all $n\gg 0$.
\end{enumerate}
For admissible modules we showed 
\cite{SchlShSug} (see also \cite{SchlSug})
that there exists a Sugawara construction which yields a 
representation of the centrally extended vector field algebra.
Here I can only give the main ideas.
Choose $(u_i)_{i=1,\ldots,\dim\g}$ and 
$(u^i)_{i=1,\ldots,\dim\g}$
a system of dual basis of $\g$ with respect to the bilinear form $(.|.)$.
The current associated to the basis element $u_i$  is given as
\begin{equation}
  J_i(Q)=\sum_{n\in\Z}\sum_{p=1}^N u_i(n,p)\w^{n,p}(Q)\ ,\quad Q\in M.
\end{equation}
Here we used $u_i(n,p)$ to denote the operator corresponding to 
$u_i\otimes A_{n,p}$ on $\Wh$.
The current is a formal operator-valued $1$-differential on  $M$.
The similar expression is used for 
the current $J^i(Q)$ associated to $u^i$, resp. for 
the current associated to an arbitrary
element $x\in\g$.
The energy-momentum tensor is defined as the formal sum
\begin{equation}
  \label{E:emtensor}
  T(Q):=\frac 12\sum_{i=1}^{\dim\g}\nord{J_i(Q)J^i(Q)}\ , \quad Q\in M\ .
\end{equation}
In this expression $\nord{...}$ denotes some normal ordering, which
moves the positive degree elements to the right.
Using the admissibility and the normal ordering we can conclude
that the energy momentum tensor is indeed a well-defined 
formal operator-valued $2$-differential on  $M$
written as 
\begin{equation}
  \label{E:emexp}
  T(Q)=\sum_{k\in\Z}\sum_{r=1}^N L(k,r)\Omega^{k,r}(Q)\ , \quad Q\in M.
\end{equation}
Let $\ka$ be the dual Coxeter number (i.e. $2\ka$ is the eigenvalue
of the Casimir operator in the adjoint representation), then 
\begin{theorem}
\cite{SchlShSug}
Assume $\g$ to be simple or abelian and $\ka$ the dual Coxeter number, resp.
$\ka=0$ for the abelian case. Assume $\Wh$ to be an admissible module
of $\gh$ with level $c$, then the
$L_{k,r}$ are well-defined operators on $\Wh$.
If $c+\ka\ne 0$ then the rescaled operators 
\begin{equation}
  L_{k,r}^*=-\frac {1}{c+\ka}\,\,L_{k,r}
\end{equation}
define a representation of a centrally extended vector
field algebra $\Lh$.
\end{theorem}
It is shown in \cite{SchlShSug} that a different normal ordering
yields an equivalent central extension.
\section{Moduli spaces and formal \\
Knizhnik-Zamo\-lod\-chi\-kov equations} 
\label{S:moduli}
The construction done in the previous sections globalizes over a
dense subset of the moduli space $\MgN$ of smooth projective curves 
of genus $g$ with 
$N$ marked points.
Recall that  two configurations 
$m=(M,(P_1,\ldots,P_N))$ and $m'=(M',(P_1',\ldots,P_N'))$
describe the same point in  moduli space if there is an 
algebraic isomorphism $\vp:M\to M'$ which maps
$P_i\to P_i'$ for $i=\iN$.
``Bad'' points in moduli space correspond to configurations 
which admit nontrivial automorphisms.
At such points singularities can occur.
For $N\ge 1$ and $g\ge 2$ the generic moduli point does not admit any
nontrivial automorphism.
Over the locus $Y$ of these points there exists a universal family of curves
with marked points
(a versal family would suffice).
The situation for $g=0$ and $g=1$ is a little bit different. Here one
should better work with the configuration space. 
This situation will not be covered here, see \cite{SchlShSug}.

Beside $\MgN$ we have also to consider $\MgNe$. 
By forgetting the last marked point
 we obtain a morphism 
$f:\MgNe\to\MgN$.
Denote by $Y$ the dense open subset of  $\MgN$, where the universal family
exists and set $\Yt=f^{-1}(Y)$ then there exists also a universal 
family over  $\Yt$.
Recall that a universal family over $\Yt$ consists of
a suitable well-behaved morphism (i.e. proper and flat)  
$\pi:\U\to \Yt$ 
for 
\newline
$\bt=[\mpt]\in\Yt$
with 
$\pi^{-1}(\bt)\cong M$, and
$N+1$ sections
\begin{equation}
  \sigma_i:\Yt\to\U,\quad \text{with}\quad
\sigma_i(\bt)=P_i\in \pi^{-1}(\bt),\qquad
i=\iNf\ .
\end{equation}
Note that for $g\ge 2$ the family $\U$ can be obtained
by pulling back the universal curve defined over an open dense subset of 
$\mathcal{M}_{g,0}$ 
via the morphism corresponding to forgetting the points.
For $g=2$ we have to start with 
$\mathcal{M}_{2,1}$.

Now we fix a reference section $\sinf$.
This corresponds to choosing a reference point for every $M$ 
depending algebraically on  varying  $M$.
Inside of $\Yt$ we have the subset
\begin{equation}
  Y'=\{[\mpt]\in\Yt\mid \Pif=\sinf(M)\}\ .
\end{equation}
By genericity the forget morphism restricted to $Y'$ 
is 1:1 onto $Y$.
We will identify in the following $Y'$ with $Y$. But note that this 
identification will not be canonical. It depends on the reference section
chosen.

Let $\pi:\U\to Y$ be the projection. 
For every open subset  $U$  of $Y$
 set 
\begin{equation}
  \St_U:=\sum_{i=1}^N\sigma_i(U)+\sinf(U)
\end{equation}
for the divisor of sections on $\uni_{|\pi^{-1}(U)}$.
The sheaf $\A_Y$ is defined as the sheaf over $Y$ given by 
defining $\A_Y(U)$ to be the space of functions on 
$\pi^{-1}(U)$ with possible poles along $\St_U$.
The same works for the central extension
$\Ah_Y$, the loop algebra $\gb_Y=\A_Y\otimes \g$ and the affine algebra
$\gh_Y=\gb_Y\oplus \mathcal{O}_Y\cdot t$.
Here $\O_Y$ denotes the structure sheaf of $Y$.
The central extensions are given in a natural globalization of
\refE{centf}, resp. \refE{knaff}.
\begin{definition}
  A sheaf $\mathcal{W}$ of $\mathcal{O}_Y$-modules is called a sheaf
of representations for the affine algebra sheaf $\gh_Y$ if 
the $\mathcal{W}(U)$ are modules over $\gh_Y(U)$ for every $U$ 
fulfilling the obvious compatibility conditions on the sheaf level.
\end{definition}
Admissible representation sheaves are defined in a similar way.
The Verma module construction can be made on the sheaf level. This yields
admissible representation sheaves  $\Wh_{(\l,\delta),Y}$.
The sheaf of conformal blocks
is the quotient sheaf
\begin{equation}
  \label{E:shcb}
  \Vh_{(\l,\delta),Y}:=
   \Wh_{(\l,\delta),Y}/\gh^*_{-,Y}\Wh_{(\l,\delta),Y}\ .
\end{equation}

Next we want to describe the tangent space at a moduli point.
\newline
Let $b=[\mpp]$ be the moduli point and set
$S=\sum_{i=1}^N P_i$, then the Kodaira-Spencer map gives an isomorphism
of the tangent space to the moduli space with certain cohomology spaces
\begin{equation}
  \label{E:KS}
  \tang_{[M]}\mathcal{M}_{g,0}\cong \mathrm{H}^1(M,T_M),\qquad
  \tang_b\mathcal{M}_{g,N}\cong \mathrm{H}^1(M,T_M(-S)) \ .
\end{equation}
Here $T_M$ is the holomorphic tangent line bundle on $M$, i.e. 
$T_M=\can^{-1}$.
The following theorem proven in \cite{SchlShWZ1}
gives an identification of the tangent space  at 
the moduli point $b$ with a certain subspace of the critical strip
of the vector field algebra associated to  the geometric data
$(M,(P_1,\ldots,P_N,\sinf(M))$.
\begin{theorem}
  There exists a natural isomorphism
\begin{equation}
(\La_{k-2}\oplus\La_{k-3}\cdots \oplus\La_{-1})\oplus \La_{(0)}^*\cong
\mathrm{H}^1(M,T_M(-kS))\ ,\quad k\ge 0\ .
\end{equation}
\end{theorem}
\noindent
Here  $\La_{(0)}^*$ are the elements of the reduced critical strip
generated by the basis elements with poles at $\PNi$.
For $g\ge 2$ its dimension is $3g-3$ and it starts with $\La_{-2}$
from above.

Let me only indicate the construction of these isomorphisms.
It is based on the calculation of the Cech cohomology
of $T_M(-kS)$ with respect to the affine (resp.
Stein) covering of $M$ given as follows.
Let $U_\infty$ be a disk around $\Pif$ and set $U_1:=M\setminus S$.
Note that $U_1$ is affine.
Then $U_1\cap U_\infty=U^*_\infty$ is the disc with $\Pif$
removed.
Hence the  Cech 1-cocycles can be given as sections on $U^*_\infty$.
In this way $f\in\La\to f_{|U^*_\infty}$ gives a linear map from the
vector field algebra to the Cech 1-cocycles, and further to 
the cohomology group $\mathrm{H}^1(M,T_M(-kS))$.
The restrictions of elements coming from  outside of the strip
given in the formulation of the theorem can by their very
definition be extended either to $U_\infty$ or to $U_1$
with the required zero-order at $S$.
Hence they are coboundaries. A closer examination shows that
the basis of the rest stays linearly independent in cohomology.
By dimension count it follows that the map is an isomorphism.

To give an example: for the element $\ e_{-1,p}\ $ (with $p=\iN$) the  
lowest order term of its expansion at  the point $P_p$ has the  
form  $\partial/\partial \xi_p$. At all other marked points it has a zero.
Under the isomorphism it corresponds
to moving the marked point $P_p$ on $M$.

Now take $\ X_k\in \tang_b\MgN\cong  \La_{(0)}^*\oplus \La_{-1}\ $
a tangent vector corresponding to an element $l_k$ of the critical 
strip.
We assume that $X_k$ operates linearly as operator $\partial_k$ 
on the space of sections of a representation sheaf $\mathcal{W}$.
For $\Phi$ a section we set
\begin{equation}
  \nabla_k\Phi:=(\partial_k+T[l_k])\Phi,\quad\text{with}\quad
T[l]:=\frac {-1}{(c+\ka)\, 2\pi \i }\int_C T(Q)l(Q)
\end{equation}
the Sugawara operator, which corresponds to $l$.

\begin{definition}
\label{D:KZe}
  The formal \KZ\ equations are defined to be the set of 
equations 
\begin{equation}
  \label{E:KZe}
  \nabla_k\Phi=0,\qquad k=1,\ldots,3g-3+N\ .
\end{equation}
\end{definition}
Note that these equations can be  expressed in terms of 
the geometric data of the curve and the points which can be moved.
In \cite{SchlShWZ1} the equations have been explicitely calculated for
genus 0 and 1.

\section*{Acknowledgments} 
First, let me point out that the new results
presented here are joint work with Oleg K. Sheinman.
I would like to thank him for fruitful cooperation 
continuing now over several years.
For the financial support of this cooperation I like to thank the
DFG and the RFBR.
It is a pleasure for me to thank J\"urgen Fuchs and 
Christoph Schweigert
for many discussions and useful hints on the subject.
\providecommand{\bysame}{\leavevmode\hbox to3em{\hrulefill}\thinspace}


\begin{thebibliography}{10}

\bibitem{Ber1}
D.~Bernard, \emph{On the {W}ess-{Z}umino-{W}itten models on the torus}, Nucl.
  Phys. B \textbf{303} (1988), 77--93.

\bibitem{Ber2}
\bysame,  \emph{On the {W}ess-{Z}umino-{W}itten models on {R}iemann
  surfaces}, Nucl. Phys. B \textbf{309} (1988), 145--174.


\bibitem{Fel}
G.~Felder, \emph{The {KZB} equations on {R}iemann surfaces}, Quantum
  symmetries/ Symetries quantiques. Proceedings of the Les Houches summer
  school (Les Houches, France, August 1 - September 8, 1995)
  (A.~Connes et. al., ed.), North-Holland, Amsterdam, 1998,
  pp.~687--725, hep-th/9609153,.

\bibitem{FeWi}
G.~Felder and Ch. Wieczerkowski, \emph{Conformal blocks on elliptic curves and
  the {K}nizhnik-{Z}amolodchikov-{B}ernard equation}, Commun. Math. Phys.
  \textbf{176} (1996), 133--161.

\bibitem{Hi}
N.~Hitchin, \emph{Flat connections and geometric quantization}, Commun. Math.
  Phys. \textbf{131} (1990), 347--380.

\bibitem{Iv}
D.~Ivanov, \emph{{K}nizhnik-{Z}amolodchikov-{B}ernard equations on {R}iemann
  surfaces}, Int. J. Mod. Phys. A \textbf{10} (1995), 2507--2536.

\bibitem{KacB}
V.K. Kac, \emph{Infinite dimensional {L}ie algebras}, Cambridge Univ. Press,
  Cambridge, 1990.

\bibitem{KniZam}
V.G. Knizhnik and A.B. Zamolodchikov, \emph{Current algebra and {W}ess-{Z}umino
  model in two dimensions}, Nucl. Phys. B \textbf{247} (1984), 83--103.

\bibitem{KN}
I.M. Krichever and S.P. Novikov, \emph{Algebras of {V}irasoro type, {R}iemann
  surfaces and structures of the theory of solitons}, Funktional Anal. i.
  Prilozhen. \textbf{21} (1987), 46.

\bibitem{Schlth}
M.~Schlichenmaier, \emph{Verallgemeinerte {K}richever - {N}ovikov {A}lgebren
  und deren {D}arstellungen}, Ph.D. thesis, Universit{\"{a}}t Mannheim, 1990.

\bibitem{Schlce}
\bysame, \emph{Central extensions and semi-infinite wedge representations of
  {K}richever-{N}ovikov algebras for more than two points}, Lett. Math. Phys.
  \textbf{20} (1991), 33--46.

\bibitem{Schlct1}
\bysame, \emph{Differential operator algebras on compact {R}iemann surfaces},
  Generalized Symmetries in Physics (Clausthal 1993, Germany) (H.-D. Doebner,
  V.K. Dobrev, and A.G. Ushveridze, eds.), World Scientific, 1994,
  pp.~425--434.

\bibitem{SchlSug}
\bysame, \emph{Sugawara operators for higher genus {R}iemann surfaces}, Rep. on
  Math. Phys. \textbf{43} (1999), 327--336, math.QA/980603.


\bibitem{SchlShSug}
M.~Schlichenmaier and O.K. Sheinman,  \emph{{S}ugawara construction and {C}asimir operators for
  {K}richever-{N}ovikov algebras}, Jour. of Math. Science \textbf{92} (1998),
  3807--3834, q-alg/9512016.

\bibitem{SchlShWZ1}
\bysame, \emph{{W}ess-{Z}umino-{W}itten-{N}ovikov theory,
  {K}nizhnik-{Z}amolodchikov equations, and {K}richever-{N}ovikov algebras,
  {I}.}, Russian Math. Surv. (Uspeki Math. Naukii). \textbf{54} (1999),
  213--250, math.QA/9812083.

\bibitem{SchlShWZ2}
\bysame,  \emph{{W}ess-{Z}umino-{W}itten-{N}ovikov
  theory, {K}nizhnik-{Z}amolodchikov equations, and {K}richever-{N}ovikov
  algebras, {II}.}, in preparation.


\bibitem{Shhwe}
O.K. Sheinman, \emph{Highest weight modules over certain quasigraded {L}ie
  algebras on elliptic curves}, Funktional Anal. i. Prilozhen. \textbf{26}
  (1992), 203--208.

\bibitem{Shhwe2}
\bysame, \emph{Affine {L}ie algebras on {R}iemann surfaces}, Funktional Anal.
  i. Prilozhen. \textbf{27} (1993), 54--62.

\bibitem{Shhwe3}
\bysame, \emph{Highest weight modules for affine {L}ie algebras on {R}iemann
  surfaces}, Funktional Anal. i. Prilozhen. \textbf{28} (1995), 44--55.

\bibitem{Tsuch}
Y.~Tsuchimoto, \emph{On the coordinate-free description of the conformal
  blocks}, J. Math. Kyoto Univ. \textbf{33} (1993), 205--243.

\bibitem{TUY}
A.~Tsuchiya, K.~Ueno, and Y.~Yamada, \emph{Conformal field theory on universal
  family of stable curves with gauge symmetries}, Adv. Stud. Pure Math.
  \textbf{19} (1989), 459--566.

\bibitem{Ucft}
K.~Ueno, \emph{Introduction to conformal field theory with gauge symmetries},
  Geometry and Physics, Proceed. Aarhus conference 1995 (Andersen~J.E. et. al.,
  ed.), Marcel Dekker, 1997, pp.~603--745.

\end{thebibliography}
\end{document}